\documentclass[a4paper,12pt]{amsart}
\usepackage{amsfonts}
\usepackage{amssymb}
\def\oC{\widehat{\mathbb{C}}}
\def\oR{\widehat{\mathbb{R}}}
\def\C{\mathbb{C}}
\def\R{\mathbb{R}}
\def\RP{\mathbb{R}\mathrm{P}}

\def\N{\mathbb{N}}
\def\veps{\varepsilon}
\def\vphi{\varphi}
\def\im{\mathrm{Im\ }}
\def\re{\mathrm{Re\ }}
\def\d{\partial}
\def\int{\mathrm{Int\ }}
\def\id{\mathrm{id}}

\newtheorem{theorem}{Theorem}[section]
\newtheorem{lemma}{Lemma}[section]
\newtheorem{corollary}{Corollary}[section]
\newtheorem{theorN}{Theorem}

\title{A compactification and the Euler characteristic of the
spaces of real meromorphic functions}

\author{S. V. Shadrin}

\thanks{partially supported by RFBR, grants 01-01-00660
and 00-15-96084}

\begin{document}

\begin{abstract}

For any connected component $H_0$ of the space of real meromorphic
function we build a compactification $N(H_0)$ of the space $H_0$.
Then we express the Euler characteristics of the spaces $H_0$ and
$N(H_0)$ in terms of topological invariants of functions from $H_0$.

\end{abstract}

\maketitle

\tableofcontents

\section{Introduction}

\subsection{Defenition of a RMF}
Recall, that a \emph{real meromorphic function (RMF)} is a triple
$(P,\tau,f)$, where $P$ is a compact Riemann surface, $\tau\colon
P\to P$ is an antiholomorphic involution, and $f\colon P\to\oC$ is a
holomorphic map to the Riemann sphere $\oC=\C\cup\{\infty\}$, s.~t.
$f(\tau(z))=\overline{(f(z))}$.

Each polynomial $p(z)=a_0+a_1z+\dots+a_nz^n$ with real coefficients
($a_i\in\R$, $i=0,\dots,n$) can be considered, for example, as a real
meromorphic function $(\oC, \eta, p)$, where $\eta\colon z\mapsto
\overline z$ is a standart involution on the Riemann sphere.

\subsection{Space of RMFs}
Denote by $H$ the set of all real meromorphic functions, i.~e. the set
of all possible triples $(P,\tau,f)$. The set $H$ posseses a natural
structure of topological space. We recall it's definition later.  The
space $H$ consists of the countable number of (arcwise) connected
components. Each component is a non-closed real manifold.

In~\cite{n1, n2, n3} one can find the full topological
classification of real meromorphic functions. There are defined
some integer topological invariants of RMFs, s.~t. each connected
component of the space $H$ is uniquely defined by assignment of these
invariants.

Consider a real meromorphic function $(P,\tau,f)$. The simplest
examples of such invariants are degree $n$ of the map $f$ and genus
$g$ of the curve $P$.

\subsection{Compactification}
For each connected component $H_0$ of the space $H$ we build a
compactification $N(H_0)$. Our construction is a real analogue
of the Natanzon-Turaev compactification~\cite{nt}.

The main idea is to represent a typical limit point of $H_0$ by a
simple branched RMF endowed with a family of disjoint disks in the
image. It is just the idea of Natanzon and Turaev, but some
refinements of the real case make our definitions a little bit
unwieldy.

The space $N(H_0)$ has a natural structure of CW-complex.
A \emph{real decorated function} is a term that we use to call an
element of the space $N(H_0)$.

Although there are many constructions of compactifications of the
spaces of non-real meromorphic functions (see~\cite{de, hm, ko}),
our work seems to be the first one that contains compactifications
of the spaces of RMFs.

\subsection{Expressions for topological characteristics}
Consider a real meromorphic function $(P,\tau,f)$. Let $H_0$ be a
connected component contains the function $(P,\tau,f)$. The space
$N(H_0)$ is a compactification of $H_0$.

Let us know all topological invariants of the function $(P,\tau,f)$.
Then there is a Riemann-Hurwitz formula for dimension:  $\dim H_0=\dim
N(H_0)=2g+2n-2$.

In this paper we express the Euler characteristics $\chi(H_0)$
and $\chi(N(H_0))$ of the spaces $H_0$ and $N(H_0)$ in topological
inveriants of the function $(P,\tau,f)$. Under the Euler
characteristic of a non-compact manifold $H_0$ we understand the
difference $\chi(H_0)=\chi(N(H_0))-\chi(N(H_0)\setminus H_0)$ (the
complement $N(H_0)\setminus H_0$ is also a CW-comlex).

\subsection{A way to compute the Euler characteristic}
We compute the Euler characteristics in two steps.

First we relate the Euler characteristics of the spaces $H_0$ and
$N(H_0)$ with the multiplicity of the \emph{Lyashko-Looijenga
map}.

Recall, that the \emph{Lyashko-Looijenga map ($LL$)}, as it was
called in~\cite{a}, associates to a real meromorphic function the
unorder sequence of it's critical values. A concept of the set of
critical values can be naturally defined also for real decorated
functions.

Then, roughly speaking, the Euler characterisic of the
manifold $H_0$ (or of the CW-complex $N(H_0)$) is equal to the number
of real meromorphic (resp., real decorated) functions from
$H_0$ (resp., from $N(H_0)$) with the set of critical values
equals to $\{i,-i\}$.

The second step is to express the obtained multiplicity as a number.

\subsection{Calculation of the Euler characteristics}
The Euler characteristic of a connected component is equal
$0$ or $1$. We explain, how it depends on topological invariants of
RMFs.

For the Euler characteristic of a compactification we don't have
a direct combinatorial formula. We express it as a number of certain
graphs.  The structure of these graphs depends on topological
invariants of RMFs from the corresponding connected component.

\subsection{Organization of the paper}
The paper is organized as it follows:

\begin{description}

\item[sect.~2] We briefly formulate the results of~\cite{n3} about
topological invariants of RMFs (as it was done in~\cite{ns} for some
other purposes).

\item[sect.~3] We define a construction of the compactification.

\item[sect.~3] We define and study the Lyashko-Looijenga map.

\item[sect.~4] We formulate and prove compactification theorems.

\item[sect.~6,~7] We represent the Euler characteristic as the
multiplicity of the Lyashko-Looijenga map.

\item[sect.~8-10] We calculate the Euler characteristics.

\end{description}

At the end of the almost all sections (except sections 5, 7 and 9)
there is a subsection contains an example. As a way of example we
always use one and the same connected component.

\subsection{Acknowledgements}
The author is grateful to S.~M.~Natanzon for numerous helpful remarks
and discussions.


\section{Real meromorphic functions}\label{rmf}

\subsection{Definition of a RMF}
As we have already said, a \emph{real meromorphic function (RMF)} is
a triple $(P,\tau,f)$, where $P$ is a compact Riemann surface,
$\tau\colon P\to P$ is an antiholomorphic involution (i.~e.,
$(P,\tau)$ is a real algebraic curve), and $f\colon P\to\oC$ is a
holomorphic map to the Riemann sphere $\oC=\C\cup\{\infty\}$, s.~t.
$f(\tau(z))=\overline{(f(z))}$.

Functions $(P_1,\tau_1,f_1)$ and $(P_2,\tau_2,f_2)$ are identified,
iff there is a biholomorphic map $\vphi\colon P_1\to P_2$, s.~t.
$f_1=f_2\vphi$ and $\vphi\tau_1=\tau_2\vphi$.

\subsection{Real algebraic curves}
Consider a real algebraic curve $(P,\tau)$. The set of real points of
the curve $(P,\tau)$ is the set $P^\tau$ of fixed points of the
involution $\tau$. It consists of pairwise disjoint simple closed
contours, called \emph{ovals}.

If $P\setminus P^\tau$ is a disconnected (connected) set, then the
curve is said to be \emph{separating} (resp., \emph{non-separating}).

The \emph{topological type of a curve $(P,\tau)$} is a triple of
numbers $(g,k,\veps)$, where $g$ is the genus of $P$, $k$ is the
number of ovals, and $\veps$ is equal to $1$ for separating curves
and to $0$ for non-separeting ones.

\subsection{Functions on non-separating curves}
Consider a function $(P,\tau,f)$. Let $c$ be an oval of the curve
$(P,\tau)$. We can choose orientations on the oval $c$ and on the
contour $\oR=\R\cup\{\infty\}\subset\oC$. Then, since
$f(c)\subset\oR$, one can define the degree of the map $f|_c\colon
c\to\oR$.

The \emph{index of the function $(P,\tau,f)$ on the oval $c$}
is an absolute value of the degree of the map $f|_c\colon
c\to\oR$. Evidently, the index doesn't depends on the choice of
orientations.

The \emph{topological type of a real meromorphic function
$(P,\tau,f)$ defined on a non-separting curve} is a tuple of numbers
$(g,n,0|I)$, where $(g,k,0)$ is the topological type of the curve
$(P,\tau)$, $n$ is the number of sheets of the covering $f$, and
$I=(i_1,\dots,i_k)$ are the indices of the function $(P,\tau,f)$ on
the ovals of the curve $(P,\tau)$.

We denote by $H(g,n,0|I)$ the set of all RMFs of topological type
$(g,n,0|I)$.

\begin{theorN}
\cite{n1, n3}
The set $H(g,n,0|I)$ is not empty, iff
$0\leq k\leq g$; $\sum^k_{j=1}i_j\leq n-2$; and
$\sum^k_{j=1}i_j\equiv n\pmod{2}$.
In this case $H(g,n,0|I)$ is a connected real manifold of dimension
$2(g+n-1)$.
\end{theorN}

\subsection{Functions on separating curves}
Consider a function $(P,\tau,f)$ defined on a separating curve.

We fix orientation on the sphere $\oC$ and on the surface $P$ as
on a complex manifold. Also we fix orientation on the contour $\oR$
as on the boundary of the upper hemisphere ($\oR=\d\Lambda$,
$\Lambda=\{z\in\oC,\ \im z>0\}$).

Then we can choose a connected component $P_1$ of the complement
$P\setminus P^\tau$. The orientation of the surface $P_1$ induces
an orientation on it's boundary $P^\tau=\partial P_1$.

Since $f(P^\tau)\subset\oR$, it is possible to define the
\emph{degrees} $I=(i_1,\dots,i_k)$ of the function $(P,\tau,f)$ on
the ovals of the real curve $(P,\tau)$.

The tuple of degrees is defined up to the simultaneous change of the
signs of all degrees, i.~e. up to the substitution $I\to -I$,
$(i_1,\dots,i_k)\to(-i_1,\dots,-i_k)$.
It is related to the arbitrariness of the choise of a connected
component of $P\setminus P^\tau$

The \emph{topological type of a real meromorphic function
$(P,\tau,f)$ defined on a separting curve} is a tuple of numbers
$(g,n,1|I)$, where $(g,k,1)$ is the topological type of the curve
$(P,\tau)$, $n$ is the number of sheets of the covering $f$, and
$I=(i_1,\dots,i_k)$ is the tuple of degrees of the function
$(P,\tau,f)$ on the ovals of the curve $(P,\tau)$.

The topological type is defined up to the substitution
$(g,n,1|I)\mapsto (g,n,1|-I)$. The topological type $(g,n,1|I)$ is
said to \emph{admit extension}, if
$|\sum^k_{j=1}i_j|<\sum^k_{j=1}|i_j|=n-2.$

We denote by $H(g,n,1|I)$ the set of all RMFs of topological type
$(g,n,1|I)$. It follows from our definitions, that
$H(g,n,1|I)=H(g,n,1|-I)$.

\begin{theorN}\cite{n2, n3}
The set $H(g,n,1|I)$ is not empty, iff
$1\leq k\leq g+1$;
$k\equiv g+1\pmod{2}$;
$\sum^k_{j=1}i_j\equiv n\pmod{2}$;
and one of the following conditions is valid:
\begin{enumerate}
\item $n=1$, $g=0$, $i_1=\pm 1$;
\item $n=2$, $k=g+1$, $i_1=\dots=i_k=0$;
\item $n\geq 2$, $|\sum^k_{j=1}i_j|=\sum^k_{j=1}|i_j|=n$, $i_j\ne 0$;
\item $n\geq 3$, $\sum^k_{j=1} |i_j|\leq n-2$.
\end{enumerate}

If the set $H(g,n,1|I)$  is non-empty and the topological type
$(g,n,1|I)$ admits no extension, then $H(g,n,1|I)$ is a connected real
manifold of dimension $2(g+n-1)$.
\end{theorN}

\subsection{Extended topological type}
Now suppose that the topological type $(g,n,1|I)$ admits extension
and $H(g,n,1|I)$ is a non-empty set.  For every function
$(P,\tau,f)\in H(g,n,1|I)$ we consider the connected component $H_0$
of $H(g,n,1|I)$ containing $(P,\tau,f)$. Then there is a function
$(P',\tau',f')\in H_0$ s.~t. $P'\cap {f'}^{-1}(\Lambda)$ is a
connected surface. Here $P'$ is a connected component of $P'\setminus
(P')^{\tau'}$ used to define the topological type $(g,n,1|I)$ and
$\Lambda\subset\oC$ is the upper hemisphere.

The extended topological type of the function $(P,\tau,f)\in
H(g,n,1|I)$ is defined as the tuple of numbers $(g,n,1|I,\xi)$, where
$\xi$ is the genus of $P'\cap {f'}^{-1}(\Lambda)$. It can be proved
that the extended topological type is well defined to within the
change $(g,n,1|I,\xi)\mapsto(g,n,1|-I,\frac{g-k+1}{2}-\xi)$.

We denote by $H(g,n,1|I,\xi)$ the subset of $H(g,n,1|I)$
consisting of all RMFs of extended topological type $(g,n,1|I,\xi)$.
It follows from our definitions, that
$H(g,n,1|I,\xi)=H(g,n,1|-I,\frac{g-k+1}{2}-\xi)$.

\begin{theorN}\cite{n2, n3}
The set $H(g,n,1|I,\xi)$ is non-empty, iff
$0\leq\xi\leq(g-k+1)/2$.
In this case $H(g,n,1|I,\xi)$  is a connected manifold of dimension
$2(g+n-1)$.
\end{theorN}

\subsection{Example}
Now we show an expamle of a connected component.

Consider real meromorphic functions
$(\oC,\eta,f_\lambda)$, $\lambda\in\Lambda$, $\Lambda=\{z\in\oC,\ \im
z>0\}$.  Here $\eta\colon z\mapsto\overline z$ is a standart
involution, and $f_\lambda$ is a two-sheets real covering of $\oC$ by
$\oC$ with two simple critical values, $\lambda$ and
$\overline\lambda$.

The set of such RMFs is a connected component of the space of real
meromorphic functions. It's easy to see, that in our notations this
set is denoted as $H(0,2,1|(2))$.

The map $(\oC,\eta,f_\lambda)\mapsto\lambda$ is a homeomorphism
$H(0,2,1|(2))\to\Lambda$. So the space $H(0,2,1|(2))$ is homeomorphic
to an open disk.


\section{Construction of compactification}\label{construct}
\subsection{Introduction}
Suppose that $H_0$ is one of connected components we discribe
in section~\ref{rmf}, i.~e. either $H_0=H(g,n,0|I)$, or
$H_0=H(g,n,1|I)$ and the topological type
$(g,n,1|I)$ admits no extension, or $H_0=H(g,n,1|I,\xi)$.

In this section we build a space $N_0$ compactifying $H_0$.

Below, in the section~\ref{prop_space_n}, we prove that $N_0$
is a compact and Hausdorff space. Also we prove there that
$N_0$ is a compactification of $H_0$.

\subsection{Real decorated functions}
A {\it real decorated function} is a quintuple
$(P,\tau,f,E,\{D_e\})$, where
\begin{enumerate}
\item $(P,\tau,f)\in H_0$ is a real meromorphic function with
$2(g+n-1)$ critical values (i.~e. all critical values of the covering
$f$ are simple).
\item $E$ is a finite subset of $\oC$, s.~t.
\begin{enumerate}
\item $\eta(E)=E$, $\eta\colon z\mapsto\overline z$,
\item the set $E$ contains no critical values of $(P,\tau,f)$.
\end{enumerate}
\item $\{D_e\}$ is a set of pairwise
disjoint closed disks, numerated by the elements $e$ of the set $E$,
s.~t.
\begin{enumerate}
\item $e\in\int D_e$,
\item each disk $D_e$ contains at least two critical
values of the function $(P,\tau,f)$,
\item there are no critical values of $(P,\tau,f)$ on the boundaries
of the disks $D_e$, $e\in E$,
\item if $\eta(e)=e'$, then $\eta(D_e)=D_{e'}$.
\end{enumerate}
\end{enumerate}

No that the set $E$ in the definition of a real decorated function
can be empty. In this case one can understand a real decorated
function as just a RMF with simple critical values.

Informally, a decorated function $(P,\tau,f,E,\{D_e\})$ should be
thought of as a limit of real meromorphic functions obtained from
$(P,\tau,f)$ by making all the critical values inside every disk
$D_e$ tend to the point $e$.

\subsection{Revision of the definition of a RMF}
We want to define some equivalence relations for real meromorphic
functions. For this purpose we need to specify our notations.

Now we consider a real meromorphic function $(P,\tau,f)$
as a covering $f\colon S_g\to\oC$ where $S_g$ is a fixed
genus $g$ two-dimensional closed real manifold.

By $P$ we understand a complex structure on $S_g$.
The involution $\tau$ is an involution of $S_g$, antiholomorphic
w.~r.~t. the complex structure $P$. The covering $f$ satisfies the
condition $f\tau=\eta f$. The complex structure $P$ is uniquely
defined by the covering $f$~\cite{ke}.

We reformulate the coinsiding condition. Two functions,
$(P_1,\tau_1,f_1)$ and $(P_2,\tau_2,f_2)$, are coinside, if there is
a homeomorphism $\vphi\colon S_g\to S_g$, s.~t.
$f_1=f_2\vphi$ and $\vphi\tau_1=\tau_2\vphi$.

\subsection{Equivalence relations}
Two real decorated functions are equivalent if one of them is
obtained from the other using a composition of the maps $A_\alpha$ and
$B_\beta$ that we describe below.

The map $A_\alpha\colon
(P_0,\tau_0,f_0,E,\{D_e\})\mapsto (P_1,\tau_1,f_1,E,\{D_e\})$
is defined by a path $\alpha\colon [0,1]\to H_0$, $\alpha\colon
t\mapsto (P_t,\tau_t,f_t)$, s.~t. on the set
$f_0^{-1}(\oC\setminus\bigcup_{e\in E}D_e)\subset S_g$
the complex structure, the involution and the covering are fixed.

The map $B_\beta\colon
(P_0,\tau_0,f_0,E,\{D_e^0\})\mapsto (P_1,\tau_1,f_1,E,\{D_e^1\})$
is defined by an isotopy $\{\beta_t\colon\oC\to\oC\}_{t\in [0,1]}$
of the identity homemorphism of the sphere $\oC$ (i.~e.,
$\beta_0=\id_{\oC}$). All $\beta_t$, $t\in[0,1]$, satisfy the
following conditions:
\begin{enumerate}
\item $\eta\beta_t=\beta_t\eta$,
\item points of the set $E$ are fixed points of $\beta_t$,
\item critical values of $(P_0,\tau_0,f_0)$ belong to
$\oC\setminus\bigcup_{e\in E}D_e^0$ are fized points of $\beta_t$,
\item $D_e^t=\beta_t(D_e^0)$ for all $e\in E$,
\item the covering $f_t\colon S_g\to\oC$ is coinside with
$\beta_tf_0$,
\item $P_t$ is just a deformation of the complex structure defined by
the deformation of the covering,
\item the involution $\tau_t\colon S_g\to S_g$ is coinside with
$\tau_0$.
\end{enumerate}

\subsection{The topology of the space of RMFs}
Now we recall the topology of the space of real meromorphic
functions defined on genus $g$ curves.

A sequence $(P_n,\tau_n,f_n)$ tends to a point $(P,\tau,f)$, if
there are a sequence $(P'_n,\tau'_n,f'_n)$ and a point
$(P',\tau',f')$, coinsiding (in the sense of our
definitions) with $(P_n,\tau_n,f_n)$ and $(P,\tau,f1)$ resp.,
s.~t. $f'_n\to f'$ as coverings $S_g\to\oC$ and
$\tau'_n\to\tau'$ as involutions of $S_g$.

\subsection{The topology of the space of real decorated functions}
Denote by $N_0$ the set of all equivalence classes of real decorated
functions. Now we define a topology of the space $N_0$.

An open neighbourhood of a point $x\in N_0$ can be defined by the
following data:
\begin{enumerate}
\item a decorated function $(P,\tau,f,E,\{D_e\})$ from the
equivalence class $x$;
\item a set of pairwise disjoint closed disks $\{B_l\subset
\oC\setminus\bigcup_{e\in E}D_e\}$ s.~t.
\begin{enumerate}
\item each disk $B_l$ contains exactly one critical value of
$(P,\tau,f)$,
\item if critical value of $(P,\tau,f)$ lies in
$\oC\setminus\bigcup_{e\in E}D_e$, then it belongs to a disk $B_l$,
\item the real condition for the disks $B_l$ is valid:
$\eta(\bigcup B_l)=\bigcup B_l$.
\end{enumerate}
\end{enumerate}
The corresponding neighbourhood of the point $x\in N_0$ consists of
equivalence classes of real decorated functions
$(P',\tau',f',E',\{D_{e'}\})$ s.~t.
\begin{enumerate}
\item $\bigcup_{e'\in E'}D_{e'}\subset\bigcup_{e\in E}D_e$;
\item the covering $f'$ can be represented as a composite $f'=\psi f_\alpha$, where
\begin{enumerate}
\item $\psi$ is a homeomorphism of the sphere $\oC$,
\item $\eta\psi=\psi\eta$,
\item $\psi$ is identity outside of $\bigcup B_l\cup\bigcup_{e\in
E}D_e$,
\item $(P_\alpha,\tau_\alpha,f_\alpha,E,\{D_e\})=
A_\alpha((P, \tau, f, E, \{D_e\}))$ for some (arbitrary) $\alpha$.
\end{enumerate}
\end{enumerate}
Here $P'$ is defined by the covering $f'$. The involutions $\tau'$
and $\tau_\alpha$ are coinside as involutions of $S_g$.

A subset $X\subset N_0$ is open, if for any point $x\in X$
there is an open neighbourhood $U(x)$ of $x$ s.t. $U(x)\subset X$.

Using the same arguments as in~\cite{nt}, one can prove that open
subsets define a topology on the space $N_0$.

\subsection{Notations}
If $H_0=H(g,n,0|I)$ ($H(g,n,1|I)$, $H(g,n,1|I,\xi)$, then we denote
$N_0$ as $N(g,n,0|I)$ (resp. $N(g,n,1|I)$, $N(g,n,1|I,\xi)$.

\subsection{Example}
A natural compactifications of the open disk
$\Lambda=\{z\in\oC,\ \im z>0\}$ is the closed disk
$'(\Lambda)=\{z\in\oC,\ \im z\geq 0\}$.

We show, that there is a homeomorphism
$N(0,2,1|(2))\to C(\Lambda)$.

Really, when $H_0=H(0,2,1|(2))$, then a real decorated function with
a non-empty set $E$ is just a function
$(\oC,\eta,f_\lambda,\{\mu\},\{D_\mu\})$, where
$\mu\in\oR$ and $\lambda\in D_\mu$.
Two such functions are equivalent, iff the point $\mu$ is the same
for both of them.

A homeomorphism $N(0,2,1|(2))\to C(\Lambda)$ is defined by
the formulas:
$$(\oC,\eta,f_\lambda,\{\mu\},\{D_\mu\})\mapsto\mu,\quad
(\oC,\eta,f_\lambda,\emptyset,\emptyset)\mapsto\lambda.$$


\section{Lyashko-Looijenga map}

\subsection{Stratification of $\RP^m$}
We denote by $\Pi^m$ the space of real tuples of $m$ unordered
not necessary distinct points of $\oC$. A tuple $(r_1,\dots,r_m)$ is
real, if $(\eta(r_1),\dots,\eta(r_m))$ is the same tuple, $\eta\colon
z\mapsto\overline z$. An example of a real tuple of four points is
$(i,i,-i,-i)$. Obviously, one can identify $\Pi^m\cong\RP^m$.

We represent $\Pi^m$ as a union of disjoint submanifolds.
Let $P=(p_1, \dots, p_r)$ and $Q=(q_1, \dots, q_s)$
be non-decreasing positive integer sequences, possibly empty,
s.~t. $m=\sum_{j=1}^rp_j+2\sum_{j=1}^sq_j$. We denote by $\Pi(P|Q)$
the set of tuples consist of $r\geq 0$ pairwise distinct points of
$\oR$, $s\geq 0$ pairwise distinct points of $\Lambda=\{z\in\oC,\
\im z>0\}$ and also $s$ points of $\eta(\Lambda)$, conjugated to
those of $\Lambda$. Here we get points of $\oR$ with multiplicities
$p_1,\dots,p_r$ and points of $\Lambda$ and $\eta(\Lambda)$ with
multiplicities $q_1,\dots,q_s$.

A set $\Pi(P|Q)$ is a submanifold of $\Pi^m$. Dimension of
$\Pi(P|Q)$ is equal to $r+2s$. $\Pi(P|Q)$ is not necessary
connected.

Obviously, $\Pi^m=\bigsqcup\Pi(P|Q)$.

\subsection{Lyashko-Looijenga map on $N_0$}
Consider a connected component $H_0$, that contains
$n$-sheeted functions on genus $g$ curves. Let $N_0$
be the corresponding space of equivalence classes of real decorated
functions

The \emph{Lyashko-Looijenga map}
$ll\colon N_0\to \Pi^{2(n+g-1)}$ associates to a
real decorated function $(P,\tau,f,E,\{D_e\})$
the tuple of points, consists of critical values of
$(P,\tau,f)$ lie in $\oC\setminus\bigcup_{e\in E}D_e$
and of points $e\in E$ with multiplicities equal
to the number of critical values of $(P,\tau,f)$ lie in $D_e$.

\begin{lemma}\label{ll_cont}
The Lyashko-Looijenga map is continuous.
\end{lemma}

\begin{lemma}\label{ll_cover}
The restriction of the Lyashko-Looijenga map to
the preimage of any connected component of any submanifold
$\Pi(P|Q)$ is a finite-sheeted non-ramified covering.
\end{lemma}

\subsection{Lyashko-Looijenga map on $H_0$}
A map $ll_H\colon H_0\to \Pi^{2(g+n-1)}$ associates
to a real meromorphic function the tuple of it's critical values
(w.~r.~t. multiplicities) is also called the \emph{Lyashko-Looijenga
map}.

\begin{lemma}
The map $ll_H$ is continuous.
\end{lemma}

\begin{lemma}
The restriction of the map $ll_H$ to the preimage of any
connected component of any submanifold $\Pi(P|Q)$ is a finite-sheeted
non-ramified covering.
\end{lemma}

\subsection{Some words about proves}
All lemmas of this section are evidently follow from our defenitions.

\subsection{Example}\label{llexample}
Let $H_0=H(0,2,1|(2))$.
Then the Lyashko-Looijenga map $ll\colon N(0,2,1|(2))\to\Pi^2$
is as follow:
$$(\oC,\eta,f_\lambda,\{\mu\},\{D_\mu\})\mapsto(\mu,\mu),\quad
(\oC,\eta,f_\lambda,\emptyset,\emptyset)\mapsto(\lambda,\overline\lambda).$$

We represent $\Pi^2$ as a union of submanifolds $\Pi(P|Q)$:
$$\Pi^2=\Pi(2|\emptyset)\sqcup\Pi(1,1|\emptyset)\sqcup\Pi(\emptyset|1).$$
In our case the map $ll$ is a one-sheeted covering over
$\Pi(\emptyset|1)$ and $\Pi(2|\emptyset)$
and a zero-sheeted covering $\Pi(1,1|\emptyset)$.


\section{Compactification theorems}\label{prop_space_n}

\subsection{Properties of $N_0$}
The following theorem is proved using the same arguments as
it's complex analogue in~\cite{nt}.

\begin{theorem}\label{space_N}
$N_0$ is a compact and a Hausdorff space.
\end{theorem}

\begin{proof}[The sketch of the proof.]

To proof that $N_0$ is a Hausdorff space we just need to build
disjoint neighbourhoods of two arbitrary points $x,y\in N_0$
s.~t. $ll(x)=ll(y)$. It's easy to do it using lemmas~\ref{ll_cont}
--~\ref{ll_cover}.

Compactness is proved by selection a convergent subsequence in an
arbitrary sequence of points of $N_0$. Using lemma~\ref{ll_cont}
and the fact that $\Pi^{2(g+n-1)}$is a compact space, we find
an $ll$-image of a limit point. Then from the lemma~\ref{ll_cover}
it follows the existence of the limit point.
\end{proof}

\subsection{Embedding $H_0\to N_0$}
Now we discribe an embedding of $H_0$ to $N_0$.
Denote by $H_0^*$ a subspace of $H_0$, consists of all functions
with simple critical values. Obviously, $H_0^*$
is an open and an everywhere dense subspace of $H_0$. The map
$C^*\colon H^*_0\to N_0$ associate to a function
$(P,\tau,f)$ the decorated function $(P,\tau,f,\emptyset,\emptyset)$.

\begin{theorem}\label{injection}
\begin{enumerate}
\item The map $C^*$ is a continuous embedding $H^*_0\to N_0$.
\item There is a continuous extension of the map $C^*$
to the map $C\colon H_0\to N_0$.
\item The map $C$ is a continuous embedding $H_0\to N_0$.
\item The image $C(H_0)$ is an open and an everywhere dense
subspace of $N_0$.
\end{enumerate}
\end{theorem}

\begin{proof}[The sketch of the proof]
From our definitions it follows that $C^*$ is a continuous map.

Now we construct the map $C\colon H_0\to N_0$ --
the continuous extension of the map $C^*$.

Let $(P_0,\tau_0,f_0)\in H_0$ be a function with non-simple critical
values. Denote by $E$ the set of it's non-simple critical values.
Evidently, $\eta(E)=E$. Let us choose pairwise disjoint closed disks
$\{D_e\}_{e\in E}$ s.~t. $e\in D_e$, $D_e$ doesn't contain
any other critical values of
$(P_0,\tau_0,f_0)$ and if $\eta(e)=e'$, then
$\eta(D_e)=D_{e'}$.

We define $'\colon (P_0,\tau_0,f_0)\mapsto
(P_1,\tau_1,f_1,E,\{D_e\})$, where
$(P_t,\tau_t,f_t)$, $t\in [0,1]$, is an arbitrary path
in the space $H_0$ s.~t. on the set
$f^{-1}(\oC\setminus\bigcup_{e\in E}D_e)\subset S_g$
complex structure, involution and covering are fixed
and all critical values of
$(P_1,\tau_1,f_1)$ are simple.

From the definition of equvalence relations of decorated functions
it follows, that the map $C$ is well defined.
It is also obvious that $C$ is a continuous map and
$C|_{H_0^*}=C^*$. The statement, that the image
$C(H_0)$ is an open and an everywhere dense subspace of $N_0$
one can obtain just from the definition of the topology of $N_0$.
\end{proof}

\subsection{Two Lyashko-Looijenga maps.}
Now we have two Lyashko-Looijenga maps on
$H_0$: the map $ll_H$ and the map $ll\circ C$.
Obviously, these maps are coinside,
$ll_H\equiv ll\circ C$.


\section{The Euler characteristic via $LL$-multiplicity -- I}

\subsection{A theorem on the Euler characteristic}

Denote by $\Pi_0$ the submanifold $\Pi(P_0|Q_0)$ of $\Pi^{2(g+n-1)}$,
where $P_0=\emptyset$, $Q_0=(g+n-1)$.
It means that $\Pi_0$ is a submanifold of tuples consist of two
conjugated points with multiplicity $g+n-1$.
Obviously, $\Pi_0$ is homeomorphic to an open disk.

Consider a connected component $H_0$ containing $n$-sheeted functions
on genus $g$ curves. Let $N_0$ be it's compactification.

\begin{theorem}\label{eulerll}

The Euler characteristic $\chi(N_0)$ of the space $N_0$
is equal to the multiplisity of the Lyashko-Looijenga map
$ll\colon N_0\to\Pi^{2(g+n-1)}$ restricted to the preimage
of $\Pi_0$.

The Euler characteristic $\chi(H_0)$ of the space $H_0$
is equal to the multiplisity of the Lyashko-Looijenga map
$ll_H\colon H_0\to\Pi^{2(g+n-1)}$ restricted to the preimage
of $\Pi_0$.

\end{theorem}

A proof of this theorem one can find in the next section.

\subsection{Some corollaries}

The following corollaries of the Theorem~\ref{eulerll} are
obvious.

\begin{corollary}\label{iidf}
Let $g+n-1>1$. Then the Euler characteristic of the space $N_0$
is equal to the number of non-equivalent decorated functions
$(P,\tau,f,E,\{D_e\})\in N_0$ s.~t.
$E=\{i,-i\}$, $D_{\pm i}=\{z||z\mp i|\leq 1/2\}$ and all critical
values of $(P,\tau,f)$ lie in $D_i\cup D_{-i}$.
\end{corollary}

The conditions $g+n-1>1$ is related with our requirement (in the
definition of a decorated function) for each disk $D_e$ to contain
at least two critical values.

\begin{corollary}\label{iiacc}
The Euler characteristic of the space $H_0$
is equal to the number of RMFs
$(P,\tau,f)\in H_0$ with only two critical values:
$i$ and $-i$.
\end{corollary}

It's obvious, that such real meromorphic functions exist only if genus
is equal to zero. But we discuss it later.

\subsection{Example}
The example is the space $H(0,2,1|(2))$.
We know that this space is homeomorphic to an open disk.
It's compactification is homeomorphic to a closed disk.
Hence, $\chi(H(0,2,1|(2)))=\chi(H(0,2,1|(2)))=1$.

The corresponding multiplicities of the Lyshko-Looijenga map we
calculate in the subsection~\ref{llexample}, and they also are equal
to $1$.


\section{The Euler characteristic via $LL$-multiplicity -- II}

In this section we prove the Theorem~\ref{eulerll}.

\subsection{Represetation of $N_0$ as a union of cells -- I}

First we represent the space $\Pi^{2(g+n-1)}$ as a union of
disjoint open balls.

Let $\Pi$ be a connected component of
$\Pi(P|Q)\subset \Pi^{2(g+n-1)}$,
$P=(p_1, \dots, p_r)$, $Q=(q_1, \dots, q_s)$.

Obviously, there is a finite-sheeted non-ramified covering
$\Pi\to\Pi(P'|Q')\subset
\Pi^{r+2s}$, $P'=(1,\dots,1)$ ($r$ times "$1$"), $Q'=(1,\dots,1)$
($s$ times "$1$").

Let us build a cell decomposition of $\Pi(P'|Q')$.
Obviously, $\Pi(P'|Q')$ is homeomorphic to
$W^r(\oR)\times W^s(\Lambda)$. Here we denote by
$W^k(X)$ the space of $k$-tuples of pairwise disjoint points of
$X$. Hence we must build only cell decompositions of
$W^r(\oR)$ and $W^s(\Lambda)$.

\subsection{Cell decompositions of $W^r(\oR)$ and $W^s(\Lambda)$}
The space $W^k(\oR)$ we represent as a union of two cells.
One of them has dimension $k$ and consists of tuples $(x_1,\dots,x_k)$
s.~t. $x_i\not=\infty$, $i=1,dots,k$.
One can consider numbers $x_1,\dots,x_k\in\R$ with relation
$x_1<\dots<x_k$ as coordinates on this cell.
The second cell has
dimension $k-1$ and consists of tuples $(x_1,\dots,x_{k-1},\infty)$.
One can consider numbers $x_1,\dots,x_{k-1}\in\R$ with relation
$x_1<\dots<x_{k-1}$ as coordinates on this cell.

To build the cell decomposition of $W^k(\Lambda)$ we need some new
notations. Let $z_1,z_2\in\Lambda$. The record
$z_1\prec z_2$ means, that $\re z_1<\re z_2$.
The record $z_1\preccurlyeq z_2$ means, that
$\re z_1=\re z_2$ and $\im z_1<\im z_2$.

Each cell consists of tuples
$(z_1,\dots,z_k)$ s.~t. coordinates on the cell are the numbers
$(\re z_1, \im z_1, \dots, \re z_k, \im z_k)$ with a relation
of the type $z_1*z_2*\dots *z_k$, where the sign
"$*$" is either $\prec$ or $\preccurlyeq$.

So, $W^k(\Lambda)$ is represented as a union of $2^{k-1}$ cells.

\subsection{Represetation of $N_0$ as a union of cells -- II}

We have the cell decompositions of the spaces $W^r(\oR)$ and
$W^s(\Lambda)$. Their product give us the cell decomposition of
$\Pi(P'|Q')$. We can lift this cell decomposition to $\Pi(P|Q)$.
We obtain the cell decomposition of $\Pi^{2(g+n-1)}$.
Then, using Lemma~\ref{ll_cover}, we represent the space $N_0$
as a union of disjoint cells.

The cell decompositions we build doesn't give us the CW-copmlex
structure on the spaces $\Pi^{2(g+n-1)}$ and $N_0$.
But to compute the Euler characteristic we need onle an arbitrary
representation of the space as a union of disjoint open balls.

An example of the cell decomposition of such space in a
particular case one can find in~\cite{ge}.

The next Lemma is evident:

\begin{lemma}
Each cell of the cell decomposition we have build
is either a subset of $C(H_0)$ or a subset of
$N_0\setminus C(H_0)$. No cell intersects both $C(H_0)$ and
$N_0\setminus C(H_0)$.
\end{lemma}

\subsection{A proof of Theorem~\ref{eulerll}}
\begin{proof}[A proof of Theorem~\ref{eulerll}]
To prove Theorem~\ref{eulerll} we just must show
(here we use Lemma~\ref{ll_cover}), that the Euler characteristics
of all connected components of all submanifolds $\Pi(P|Q)$, except
$\Pi_0$, are equal to zero.

By the Euler characteristic we understand
here the alternated sum of the numbers of cells we have build.

Hence we need only to prove that the Euler characteristic of any
submanifold $\Pi(P'|Q')$, except $\Pi_0\subset\Pi^2$, is equal to
zero. Really, if $P'\not=\emptyset$ (i.~e., $r>0$), then
$\chi(\Pi(P'|Q'))=0$, since $\chi(W^k(\oR))=0$.
So we consider
the case $r=0$ and $s\geq 2$. Then $\Pi(P'|Q')=W^s(\Lambda)$.
There are $2^{s-1}$ cells in the cell decomposition of
$W^s(\Lambda)$. The dimension of a cell one can compute using formula
$\dim=2+2(\prec)+(\preccurlyeq)$, where $(\prec)$ and
$(\preccurlyeq)$ are the numbers of corresponding signs in the
relation that defines the cell. So, the Euler characteristic
$\chi(W^s(\Lambda))$ is equal to
$\binom{c-1}{0}-\binom{c-1}{1}+
\binom{c-1}{2}-\dots\pm\binom{c-1}{c-1}=(1-1)^{c-1}=0$.

From the foregoing we obtain our theorem, since $\chi(\Pi_0)=1$.
\end{proof}


\section{The Euler characteristic of connected components}

\subsection{Theorems}
Denote by $H_1$ a connected component $H(g,n,0|I)$.

\begin{theorem}\label{acc1}
If $g=0$, then $\chi(H_1)=1$.
If $g\not= 0$, then $\chi(H_1)=0$.
\end{theorem}

Consider a connected component $H_2=H(g,n,1|I)$
(the topological type $(g,n,1|I)$, $I=(i_1,\dots,i_k)$,
doesn't admit extension).
Recall, that if $g=0$, then $k=1$ and $I=(i_1)$.

\begin{theorem}\label{separacc}
If $g=0$ and $|i_1|=n$, then $\chi(H_2)=1$.
In any other case $\chi(H_2)=0$.
\end{theorem}

Denote by $H_3$ a connected component $H(g,k,1|I,\xi)$.

\begin{theorem}\label{acc3}
Independently of topological invariants $\chi(H_3)=0$.
\end{theorem}

\begin{proof}[Proofs of Theorems~\ref{acc1}-\ref{acc3}]
Consider a connected component $H_0$ of the space of real
meromorphic functions.

From the construction of the map $K$ (Theorem~\ref{injection})
and Corollary~\ref{iiacc} it follows, that the Euler characteristic of
the space $H_0$ is equal to the number of functions $(P,\tau,f)\in
H_0$ with only two critical values, $i$ and $-i$.

Consider such function. Since it has only one critical value in the
disk $\Lambda=\{z\in\C|\im z>0\}$, any connected component of the
preimage of $\Lambda$ is a disk. The same is true for
$\eta(\Lambda)=\{z\in\C|\im z<0\}$. Hence, the genus of the surface
$P$ is equal to zero. Also it means, that the preimage of the contour
$\oR$ is connected.

Obviously, each connected component $H_0$ contains at least one
such function. $H_0$ contains exactly one such function, iff
$g=0$, and either $\veps=0$, or $|i_1|=n$. From this
our Theorems follow. \end{proof}

\subsection{Example}
We use Theorem~\ref{separacc} to study our usual example,
the space $H(0,2,1|(2))$. We obtain, that $\chi(H(0,2,1|(2)))=1$.


\section{The Euler characteristic in non-separating case}

Denote by $N_1$ the compactification $N(g,n,0|I)$,
$I=(i_1,\dots,i_k)$.  To express the Euler characteristic
$\chi(N_1)$ of the space $N_1$ in general case we need
\emph{non-sparating graphs}.
In the definition of a non-separating graph we suppose
that all indices $i_j$ are non-zero (or there are no of them,
i.~e. $k=0$, $I=()$).

A \emph{non-separating graph} consists of the following data:
\begin{enumerate}
\item A bipartited graph with signs on vertexes and edges:
\begin{enumerate}
\item The set of vertexes $V=V_w\sqcup V_b$. This set is a
disjoint union of two subsets with the same cardinalities
($|V_w|=|V_b|$). Vertexes from $V_w$ (from $V_b$) we suppose to be
painted in white (resp., black) color.
\item The set of edges $E$.
Each edge connects vertexes of different colors.
The graph $(V,E)$ is connected.
\item A function $\zeta_V\colon V\to \N\cup\{0\}$.
We require the following:
$g=k+|E|-|V|+1+\sum_{v\in V} \zeta_V(v)$.
\item A function $\zeta_E\colon E\to \N$.
The following must be valid:
$n=\sum_{e\in E} \zeta_E(e)-\sum_{j=1}^k i_j$.
\end{enumerate}
\item A structure of root vertexes and edges:

There are fixed subsets of root vertexes
$V_w^r\subset V_w$, $V_b^r\subset V_b$.
There are exactly $k$ root vertexes of each color.
The function $\zeta_V$ vanishes on all root vertexes.

Each root vertex has exactly one edge outgoing from it.
We denote by $E_w^r$ ($E_b^r$) the set of edges outgoing of
root vertexes of the white (resp., black) color.
The function $\zeta_E$ restricted to $E_w^r$ and $E_b^r$ give us
isomorphisms $E_w^r\to I$ and $E_b^r\to I$ of the sets $E_w^r$ and
$E_b^r$ and the tuple of indices $I=(i_1,\dots,i_k)$.

\item An automotphism $\gamma$ of the graph $(V,E)$:
\begin{enumerate}
\item The automorphism $\gamma$ change colors of all vertexes.
\item $\zeta_V\circ\gamma=\zeta_V$, $\zeta_E\circ\gamma=\zeta_E$.
\item On the edges that are fixed w.~r.~t. the automorphism $\gamma$,
the function $\zeta_E$ gets even values.
\item The automorphism $\gamma$ sends root vertexes to root
vertexes.
\end{enumerate}

\end{enumerate}

As a way of example to this definiton we discribe the unique
non-separating graph for the space $N(1,3,0|(1))$:
$V=\{v_{-2}, v_{-1}, v_1, v_2\}$, $V_w=\{v_{-2},v_{1}\}$,
$V_b=\{v_{-1},v_{2}\}$, $V^r_w=\{v_{-2}\}$, $V^r_b=\{v_2\}$,
$E=\{(v_{-2},v_{-1}), (v_{-1},v_{1}), (v_1,v_2)\}$,
$E_w^r=\{(v_{-2},v_{-1})\}$, $E_b^r=\{(v_1,v_2)\}$,
$\zeta_V\equiv 0$, $\zeta_E(E^r_{w,b})=1$, $\zeta_E((v_{-1},v_{1}))=2$,
$\gamma(v_i)=v_{-i}$.

\begin{theorem}\label{nonsepar}
If $k>0$ and at least one of the indices
$i_j$ is equal to zero, then $\chi(N_1)=0$.
In any othe case the Euler characteristic
$\chi(N_1)$ is equal to the number of different
non-separating graphs.
\end{theorem}

\begin{proof}
Consider a non-empty space $H_1=H(g,k,0|I)$

Let the topological type $(g,k,0|I)$ contain zero index.
Then for any decorated function
$(P,\tau,f,E,\{D_e\})\subset N_1=N(H_1)$ we obtain the following:
either function $(P,\tau,f)$ has a real critical value in
$\oC\setminus\bigcup D_e$, or the set $E$ contains a real point,
or both statements are true. Hence, the image of the
Lyashko-Looijenga map doesn't contains the submanifold
$\Pi_0$. Then, using Theorem~\ref{eulerll}, we conclude that
$\chi(N_1)=0$.

Now we suppose that all indices $i_j$ of the topological type
$(g,k,0|I)$ are non-zero (or there are no of them, i.~e. $k=0$,
$I=()$). Then the space $N_1$ contains functions that described
in Corollary~\ref{iidf}.

To each decorated function $(P,\tau,f,E,\{D_e\})$ s.~t.
$E=\{i,-i\}$, $D_{\pm i}=\{z||z\mp i|\leq 1/2\}$ and all critical
values of $(P,\tau,f)$ are in the disks $D_{\pm i}$,
we associate a non-separating graph.

We choose an arbitrary function $(P,\tau,f)$ that defines
an equvalence class of $(P,\tau,f,E,\{D_e\})$.
We cut the surface $P$ along ovals. Then we glue all holes we obtain
using disk coverings of the disks $\Lambda$ and
$\eta(\Lambda)$, s.~t.  the result is a real meromorphic function
$(P',\tau',f')$.  Obviously, the topological type of $(P',\tau',f')$
is equal to $(g-k,n+\sum_{j=1}^k i_j,0|I')$, $I'=()$.

To each connected component of $f^{-1}(\Lambda)$
($f^{-1}(\eta(\Lambda))$) we associate a white (resp., black) vertex
with the value of $\zeta_V$ equals to the genus of this component.
The root vertexes are obtaind from the glueing disks.
To each component of $f^{-1}(\oR)$ we associate an edge
with the value of $\zeta_E$ equals to the number of sheetes of
the covering of the contour $\oR$ by this component. An edge outgoes
from a vertex if the corresponding components of the surface
$P'$ are bounded.

The automorphism of the graph is generated by the involution
$\tau'$. It's easy to see that the graph we build is non-separating.

\begin{lemma}\label{necongdf}
We associate to non-equivalent decorated functions different
non-separating graphs. Each non-separating graph can be
obtained as the graph of a function.
\end{lemma}

\begin{proof}
The second statement of the Lemma one can prove by direct
construction. The first statement of the Lemma is an obvious
corollary of Lemma~\ref{cover} which you can see below.

See also~\cite{nsv}. There are the standart arguments that allow
one to prove this Lemma are written in details.
\end{proof}

Consider two $n$-sheeted coverings $f_1,f_2\colon
U\to D$ of a disk $D$ by a connected surface $U$ with only simple
critical values. Let $f_1$ and $f_2$ be topologically equvalent over
the boundary of $D$.

\begin{lemma}\label{cover}
There are homeomorphisms $\vphi\colon U\to U$
and $\psi\colon D\to D$
s.~t. $\psi f_1=f_2\vphi$.
\end{lemma}

This Lemma is a corollary of results~\cite{n1, n3, kz} on the
topological classification of coverings of the sphere $\oC$
with one non-simple critical value.

Let us finish the proof of the Theorem. From Lemma~\ref{necongdf}
it follows, that the number of functions discribed in Corollary~\ref{iidf}
is equal to the number of non-separating graphs.
Hence the Euler characteristic of the space $N(g,k,0|I)$ is equal to
the number of non-separating graphs.
\end{proof}


\section{The Euler characteristic in separating case}

\subsection{The topological type doesn't admit extension}
Denote by $N_2$ a compactification $N(g,n,1|I)$
(the topological type $(g,n,1|I)$, $I=(i_1,\dots,i_k)$,
doesn't admit extension).
To express the Euler characteristic $\chi(N_2)$ of the
space $N_2$ in general case we need \emph{separating graphs}.
In the definition of the separating graph we suppose
that all degrees $i_j$ are non-zero.

A \emph{separating graph} consists of the following data:
\begin{enumerate}

\item A bipartite graph with signs on vertexes and edges:

\begin{enumerate}
\item The set of vertexes $V=V_w\sqcup V_b$. This set is a disjoint
union of two subsets (not necessary with the same cardinalities).
Vertexes from $V_w$ (from $V_b$) we suppose to be
painted in white (resp., black) color.

\item The set of edges $E$. Each edge connects vertexes of different
colors. The graph $(V,E)$ is connected.

\item The function $\zeta_V\colon V\to \N\cup\{0\}$.
We require the following:
$g=(k-1)+2(|E|-|V|+1)+2\sum_{v\in V} \zeta_V(v)$.
\item The function $\zeta_E\colon E\to \N$. The following must be
valid: $n=2\sum_{e\in E} \zeta_E(e)-\sum_{j=1}^k |i_j|$.
\end{enumerate}

\item The structure of root vertexes and edges:

There are fixed subsets of root vertexes
$V_w^r\subset V_w$, $V_b^r\subset V_b$.
We have exactly $k$ root vertexes at all.
The number of root vertexes of white (black) color is
equal to the number of negative (resp., positive) degrees $i_j$ in
the tuple $I=(i_1,\dots,i_k)$.
The functions $\zeta_V$ vanishes on all root vertexes.

Each root vertex has exactly one edge outgoing from it.
We denote by $E^r$ the set of edges outgoing from root vertexes.
We define the function $\zeta_R$ on the set $E_r$:
in an edge $e\in E^r$ outgoes from the root vertex of white (black)
color, then we put $\zeta_R(e)=-\zeta_E(e)$
(resp., $\zeta_R(e)=\zeta_E(e)$).  The function $\zeta_R$ must be
an isomorphism of the set $E^r$ and the tuple of degrees $I=(i_1,
..., i_k)$.

\end{enumerate}

As a way of example to this definition we discribe the unique
separating graph for the space $N(1,3,1|(1,2))$:
$V=\{v_0, v_1, v_2\}$, $V_w=\{v_1\}$, $V^r=V_b=\{v_0,v_2\}$,
$E=E^r=\{(v_0,v_1), (v_1,v_2)\}$, $\zeta_V\equiv 0$.

\begin{theorem}\label{separ}
If $|\sum^k_{j=1} i_j|=n$, then $\chi(N_2)=1$.
If at least one of the degrees $i_j$ is equal to zero, then
$\chi(N_2)=0$.  If $\sum^k_{j=1} |i_j|\leq n-2$ and $i_j\not= 0$,
then the Euler characteristic of the space $\chi(N_2)$
is equal to the number of separating graphs.
\end{theorem}

\begin{proof}
If one of the degrees $i_j$ is equal to zero, then
$\chi(N_2)=0$ by the same arguments as in the proof of
Theorem~\ref{nonsepar}.

Let all degrees $i_j$ be non-zero.
Then we associate a separating graph to each decorated function
$(P,\tau,f,E,\{D_e\})\in N_2$ s.~t. $E=\{i,-i\}$, $D_{\pm
i}=\{z||z\mp i|\leq 1/2\}$ and all critical values of $(P,\tau,f)$
are in the disks $D_{\pm i}$.

We choose an arbitrary function $(P,\tau,f)$ that defines
an equvalence class of $(P,\tau,f,E,\{D_e\})$.
Let $P_1$ be the connected component of the complement
$P\setminus P^\tau$ that we used to
define the topological type $(g,k,1|I)$.
We cut out along ovals the other part of the curve $P$.
Then we glue all holes in $P_1$ using disk coverings of
$\Lambda$ and $\eta(\Lambda)$. Thus we obtain a meromorphic function,
now not a real one. Denote by $P_1'$ the surface that we obtain.

To each connected component of $f^{-1}(\Lambda)\cap P_1$
($f^{-1}(\eta(\Lambda))\cap P_1$) we associate a white (black)
vertex with the value of $\zeta_V$ equals to the genus of the
component. Root vertexes we obtain from the glueing disks. To each
component of $f^{-1}(\oR)\cap P_1$ we associate an edge with the
value of $\zeta_E$ equals to the number of sheets of the covering
of the contour $\oR$ by this component. An edge outgoes from the
vertex, iff the corresponding components of the surface
$P_1'$ are bounded. Root edges and the isomorphism $\zeta_R$
we obtain from the tracks of oval on $P_1'$.

It's easy to see that the graph we have build is separating.

\begin{lemma}\label{separlemma}
To non-equvalent decorated functions we associate non-equivalent
separating graphs. Each separating graph can be obtained as the graph
of a function.
\end{lemma}

\begin{proof}
The first statement is a corollary of Lemma~\ref{cover}.
The second statement one can prove by direct constructon.
\end{proof}

Now we prove the third statement of the Theorem.
From Lemma~\ref{separlemma} it follow, that the number of functions
discribed in Corollary~\ref{iidf} is equal to the number of
separating graphs. Hence the Euler characteristic of the space
$N_2$ is equal to the number of
separating graphs.

If $|\sum^k_{j=1} i_j|=n$, then it's easy to see that there is a
unique separating graph.
\end{proof}

\subsection{The topological type admits extension}

Denote by $N_3$ the compactification $N(g,n,1|I,\xi)$,
$I=(i_1,\dots,i_k)$.

\begin{theorem}
If one of the degrees $i_j$ is equal to zero, then $\chi(N_3)=0$.
In any other case $\chi(N_3)=1$.
\end{theorem}

\begin{proof}
If one of the degrees $i_j$ is equal to zero, then
$\chi(N_3)=0$ by the same arguments as in the proof of
Theorem~\ref{nonsepar}. In any other case we just must show, that
there is a unique decorated function in $N_3=N(g,k,1|I,\xi)$
satisfying conditions of Corollary~\ref{iidf}.

Consider a function $(P,\tau,f)\in H(g,k,1|I,\xi)$. Suppose that
it has no critical values on the contour
$\oR$. Let $P_1$ be the component of the complement
$P\setminus P^\tau$, that we used to define the topological type
$(g,k,1|I,\xi)$.

\begin{lemma}\label{connectness}
$f^{-1}(\Lambda)\bigcap P_1$ and
$f^{-1}(\eta(\Lambda))\bigcap P_1$ are connected surfaces.
\end{lemma}

Is lemma is proved in~\cite{n3} as a little bit different statement.

From Lemmas~\ref{cover} and~\ref{connectness} it follows,
that there is a unique function in $N_3=N(g,k,1|I,\xi)$
satisfying conditions of Corollary~\ref{iidf}.
\end{proof}

\subsection{Example}
We apply Theorem~\ref{separ} to our usual example, the space
$N(0,2,1|(2))$. It's easy to see, that the separating graph is
exactly one:
$V=\{v_0, v_1\}$, $V_w=\{v_1\}$, $V^r=V_b=\{v_0\}$,
$E=E^r=\{(v_0,v_1)\}$, $\zeta_V\equiv 0$.

\bigskip \footnotesize \noindent \sc
Independent University of Moscow \\
Bolshoi Vlasjevskii pereulok 11, \\
Moscow 121002, Russia

\smallskip \noindent
e-mail: shadrin@mccme.ru


\begin{thebibliography}{00}

\bibitem{a} V.~I.~Arnold, Topological classification of trigonometric
polynomials and combinatorics of graphs with an equal number of
vertices and edges, Funct Anal. Appl. 30 (1996), 1-14

\bibitem{de} S.~Diaz, D.~Edidin, Towards the homology of Hurwitz
space, J.~Diff.~Geom.~43 (1996), 66-98

\bibitem{ge} A.~Gabrielov, A.~Eremenko, Rational functions with
real critical points and the B.~and M.~Shapiro conjecture in real
enumerative geometry, to appear in Annals of Math.

\bibitem{hm} J.~Harris, D.~Mumford, On the Kodaira dimension of the
moduli space of curves, Invent.~Math.~67 (1982), 23-86

\bibitem{ke} B.~Kerekjarto, Vorlesungen \"uber Topologie I.
Fl\"achentopologie, Springer-Verlag, Berlin, 1923

\bibitem{kz} A.~G.~Khovanskii, S.~Zdravkovska, Branched covers of
$S^2$ and braid groups, J.~Knot Theory Ramific. 5:1
(1996), 55-75

\bibitem{ko} M.~Kontsevich, Enumeration of rational curves via torus
action, in: The moduli space of curves, PM 129, Birkh\"auser, Basel,
1995, 335-368

\bibitem{n1} S.~M.~Natanzon, Spaces of real meromorphic functions on
real algebraic curves, Soviet Math. Dokl. 30 (1984), 724-726

\bibitem{n2} S.~M.~Natanzon, Real meromorphic fucntions on real
algebraic curves, Soviet Math. Dokl. 36 (1987), 425-427

\bibitem{n3} S.~M.~Natanzon,
Topology of 2--dimensional coverings and meromorphic functions on
real and complex algebraic curves,
Selecta Math. Soviet. 12:3 (1993), 251-291

\bibitem{nsv} S.~Natanzon, B.~Shapiro, A.~Vainshtein,
Topological classification of generic real rational functions,
to appear in J.~Knot Theory Ramific., math.AG/ 0110235

\bibitem{ns} S.~M.~Natanzon, S.~V.~Shadrin, Topological
classification of unitary functions of arbitrary genus, Russian Math.
Surveys 55:6 (2000), 1163-1164

\bibitem{nt} S.~Natanzon, V.~Turaev, A compactification of the
Hurwitz space, Topology 38:4 (1999), 889-914


\end{thebibliography}
\end{document}